\documentclass[a4,amstex,12pt]{article}
\usepackage{amsmath,wrapfig}
\usepackage{amsxtra}
\usepackage{amssymb,bm,epic}
\usepackage[mathscr]{eucal}
\usepackage{amsthm, amsfonts, latexsym,enumerate}
\usepackage{graphicx}
\theoremstyle{definition}
\newtheorem{dfn}{Definition}[section]
\newtheorem{thm}[dfn]{Theorem}
\newtheorem{prp}[dfn]{Proposition}

\newtheorem{crl}[dfn]{Corollary}

\newtheorem*{prf}{Proof}
\newtheorem{alg}[dfn]{Algorithm}
\theoremstyle{break}
\begin{document}
\centerline{\Large Two characteristic polynomials corresponding to}\par\vspace*{1mm}
\centerline{\Large graphical networks over min-plus algebra}\par
\vspace{4mm}
\centerline{Sennosuke WATANABE$^a$, Yuto TOZUKA$^b$, 
Yoshihide WATANABE$^c$,}\par\vspace{1mm}
\centerline{Aito YASUDA$^d$, 
Masashi IWASAKI$^e$}\par\vspace{4mm}
\centerline{{\small 
${}^a$\textit{
Department of General Education, National Institute of Technology, }}}\par
\centerline{{\small \textit{Oyama College, 771 Nakakuki, Oyama City, Tochigi, 323-0806 Japan}}}
\par\vspace{2mm}
\centerline{{\small $^b$\textit{Graduate School of Science and 
Engineering, Science of Environment and}}}
\par\centerline{{\small \textit{Mathematical Modeling, 
Doshisha University, 1-3 Tatara Miyakodani,}}}
\par\centerline{{\small \textit{Kyotanabe, 610-0394 Japan}}}
\par\vspace{2mm}
\centerline{{\small $^c$\textit{Faculty of Science and Engineering, 
Department of Mathematical Sciences,}}}\par
\centerline{{\small \textit{Doshisha University, 1-3 Tatara Miyakodani, 
Kyotanabe, 610-0394 Japan}}}\par\vspace{2mm}
\centerline{{\small ${}^d{}^e$\textit{
Faculity of Life and Environmental Sciences, 
Kyoto Prefectural University,}}}\par
\centerline{{\small \textit{
1-5 Nakaragi-cho, Shimogamo, Sakyo-ku, 
Kyoto, 606-8522 Japan}}}
\par\vspace{2mm}
\centerline{{\small\textit{E-mail addresses:} 
$^a$sewatana@oyama-ct.ac.jp, $^c$yowatana@mail.doshisha.ac.jp,}}
\centerline{{\small\hspace*{-113pt}
$^d$imasa@kpu.ac.jp}}
\begin{abstract}
In this paper, we investigate characteristic polynomials of matrices in min-plus algebra.
Eigenvalues of min-plus matrices are known to be the minimum roots of the characteristic polynomials 
based on tropical determinants which are designed from emulating standard determinants.
Moreover, minimum roots of characteristic polynomials have a close relationship to 
graphs associated with min-plus matrices consisting of vertices and directed edges with weights.
The literature has yet to focus on the other roots of min-plus characteristic polynomials.
Thus, here we consider how to relate the $2$nd, $3$rd, $\dots$ minimum roots 
of min-plus characteristic polynomials to graphical features.
We then define new characteristic polynomials of min-plus matrices 
by considering an analogue of the Faddeev-LeVerrier algorithm
that generates the characteristic polynomials of linear matrices.
We conclusively show that minimum roots of the proposed characteristic polynomials coincide with min-plus eigenvalues,
and observe the other roots as in the study of the already known characteristic polynomials.
We also give an example to illustrate the difference between the already known and proposed characteristic polynomials.
\end{abstract}
\section{Introduction}
Various fields of mathematics consider min-plus algebra, which is an abstract algebras with idempotent semirings.
The arithmetic operations of min-plus algebra are $\min(a,b)$ and $a+b$ for $a, b\in\mathbb{R}_{\min}:=\mathbb{R}\cup \{\infty\}$ 
where $\mathbb{R}$ is the set of all real numbers. 
Although it has different operators to the well-known linear algebra,
the eigenvalue problem is fundamental both types of algebra. 
The min-plus eigenvalue problem was shown in Gondran-Minoux \cite{BCOQ} and Zimmermann \cite{GM} to have a close relationship
with the shortest path problem on graphs consisting of vertexes and edges, where every edge links to two distinct vertices. 
Directions are added to edges in directed graphs, and a value is associated with each edge in weighted graphs. 
Matrices whose entries are min-plus algebra figures are sometimes considered with respect to directed and weighted graphs. 
Such matrices are called min-plus matrices, and are practically defined by 
assigning the weights of edges from the vertices $i$ to $j$ to the $(i,j)$ entries. 
According to Gondran-Minoux \cite{BCOQ} and Zimmermann \cite{GM}, if a min-plus matrix has an eigenvalue,
this eigenvalue reflects a significant feature in the network on a directed and weighted graph that is associated with the min-plus matrix.
There exist circuits whose average weights are the eigenvalue, where a circuit signifies a closed path without crossing; 
its average weight is given by the ratio of the sum of all weights to the vertex number. 
Conversely, in the network involving circuits, the minimum of the average weights of circuits 
coincides with the eigenvalue of the corresponding min-plus matrix.
Moreover, the eigenvalues of min-plus matrices are the minimum roots of the characteristic polynomials,
defined using tropical determinants, which correspond to the determinants over linear algebra, over min-plus algebra \cite{MS}.
However, the $2$nd, $3$rd, $\dots$ minimum roots have not yet been related to graphical features.
Thus, the first goal of this paper is to identify graphical significance of the $2$nd, $3$rd, $\dots$ minimum roots. 
\par
Over linear algebra, the $QR$, qd and Jacobi algorithms are representative numerical solvers for eigenvalue problem \cite{D,GL,R}. 
The divided-and-conquer and bisection algorithms are also the famous linear eigenvalue solvers \cite{D,GL}. 
In contrast, few procedures for min-plus eigenvalues have been studied, with the exception of work by Maclagan-Sturmfels \cite{MS}. 
Moreover, to the best of our knowledge, 
min-plus eigenvalue procedures based on linear equivalents have not been yet addressed in the literature.
Thus, the second goal of this paper is to propose new eigenvalue algorithms for min-plus algebra
by emulating the Faddeev-LeVerrier algorithm \cite{Fad} in linear algebra.
The Faddeev-LeVerrier algorithm employs only linear scalar and matrix arithmetic, which can be intuitively replaced with min-plus one.
Strictly speaking, the Faddeev-LeVerrier algorithm generates not eigenvalues but characteristic polynomials of square matrices. 
In other words, our second goal is to essentially show how to derive new characteristic polynomials of min-plus matrices. 
\par
The remainder of this paper is organized as follows. 
Section 2 describes elementary scalar and matrix arithmetic over the min-plus algebra.
Section 3 and 4 explain the relationships between min-plus matrices and the corresponding networks, 
and linear factorizations of min-plus polynomials, including an effective preconditioning, respectively. 
In Section 5, we elucidate not only minimum roots, but also the other roots of given min-plus characteristic polynomials 
given using tropical determinants from the perspective of graphical networks.
In Section 6, by considering an analogue of the Faddeev-LeVerrier algorithm, 
we derive new characteristic polynomials of min-plus matrices,
then clarify their features in the comparison with known characteristic polynomials. 
Finally, in Section 7, we provide concluding remarks.
%
\section{Min-plus arithmetic}
%
In this section, we present elementary definitions and properties concerning min-plus algebra.
We first focus on scalar arithmetic over the min-plus algebra, and then present the matrix arithmetic.
\par
For $a,b\in\mathbb{R}_{\min}$, min-plus algebra has only two binary arithmetic operators, $\oplus$ and $\otimes$,
which have the following definitions,
\begin{align*} 
& a\oplus b=\min\{a,b\},\\
& a\otimes b= a+b.
\end{align*}
We can easily check that both $\oplus$ and $\otimes$ are associative and commutative, 
and $\otimes$ is distributive with respect to $\oplus$, namely, for $a,b,c\in\mathbb{R}_{\min}$, 
\begin{align*}
a\otimes (b\oplus c)=(a\otimes b)\oplus (a\otimes c).
\end{align*} 
Moreover, we may regard $\varepsilon =+\infty$ and $e=0$ as identities with respect to 
$\oplus$ and $\otimes$, respectively because for any $a\in\mathbb{R}_{\min}$, 
\begin{align*}
& a\oplus\varepsilon =\min\{a,+\infty\}=a,\\
& a\otimes e=a+0=a.
\end{align*}
Using the identity $e$, we can uniquely define the inverse of $a\in\mathbb{R}_{\min}\!\setminus\!\{\varepsilon\}$ 
with respect to $\otimes$, denoted by $b$, as
\begin{align*}
a\otimes b=e.
\end{align*}
Since it holds that
\begin{align*}
a\otimes\varepsilon =a+\infty=\varepsilon,
\end{align*}
the identity $\varepsilon =+\infty$ with respect to $\oplus$ is absorbing for $\otimes$.
Here, we consider matrices whose entries are all $\mathbb{R}_{\min}$ numbers to be min-plus matrices.
For positive integers $m$ and $n$, we designate the set of all $m$-by-$n$ min-plus matrices as $\mathbb{R}_{\min}^{m\times n}$.
Since the min-plus matrices appearing in the later sections are all square matrices, 
we hereinafter limit discussion to $n$-by-$n$ min-plus matrices.
For $A=(a_{ij}),B=(b_{ij})\in\mathbb{R}_{\min}^{n\times n}$, the sum $A\oplus B =([A\oplus B]_{ij})\in\mathbb{R}_{\min}^{n\times n}$ 
and product $A\otimes B=([A\otimes B]_{ij})\in\mathbb{R}_{\min}^{n\times n}$ are respectively given as:
\begin{align*}
[A\oplus B]_{ij}=a_{ij}\oplus b_{ij}=\min\{a_{ij},b_{ij}\},
\end{align*}
and
\begin{align*}
[A\otimes B]_{ij}=\bigoplus^k_{\ell=1}(a_{i\ell}\otimes b_{\ell j})=\underset{\ell=1,2,\dots,k}{\min}\{a_{i\ell}+b_{\ell j}\}.
\end{align*}
Moreover, for $\alpha\in\mathbb{R}_{\min}$ and $A=(a_{ij})\in\mathbb{R}_{\min}^{n\times n}$, 
the scalar multiplication $\alpha\otimes A =([\alpha\otimes A]_{ij})\in\mathbb{R}_{\min}^{n\times n}$ is defined as 
\begin{align*}
[\alpha\otimes A]_{ij}=\alpha\otimes a_{ij}.
\end{align*}
%
\section{Graphs and min-plus eigenvalues}
%
In this section, we first give a short explanation for min-plus matrices corresponding to graphs which are not functional graphs.
Then, we review the relationships between the eigenvalues of min-plus matrices and the corresponding graphs.
\par
Let $v_1,v_2,\dots, v_m$ denote vertices on the graph $G$, 
and let $e_{i,j}=(v_i,v_j)$ be edges which link the vertices $v_i$ and $v_j$.
The edge $e_{i,i}=(v_i,v_i)$ is called a loop.
Moreover, let $V:=\{v_1,v_2,\dots,v_m\}$ and $E:=\{e_{i,j} | (i,j)\in\sigma\}$ 
where $\sigma$ is the set of all pairs of $i$ and $j$ such that the edge $e_{i,j}$ exists.
Then, two sets $V$ and $E$ uniquely determine the graph $G$.
Thus, such $G$ is often expressed as $G=(V,E)$.
If $G$ is a directed graph, then $e_{i,j}$ are directed edges whose tail and head vertices are $v_i$ and $v_j$, respectively.
Further, if $G$ is a directed and weighted graph, 
then the real number $w(e_{i,j})$ is assigned to each edge $e_{i,j}$, and is called the weight.
The pair $\mathcal{N}=(G,w)$ is often called the network on the graph $G$.
The following definition gives the so-called weighted adjacency matrices associated with networks.
\begin{dfn}
For the network $\mathcal{N}$ involving $m$ vertices, 
an $m$-by-$m$ weighted adjacency matrix $A(\mathcal{N})=(a_{ij})$ is given using $\mathbb{R}_{\min}$ numbers as 
\begin{align*}
a_{ij}=\left\{\begin{array}{ll}
w((v_i,v_j)) & \text{if}\ (v_i,v_j)\in E,\\
+\varepsilon & \text{otherwise} . 
\end{array}\right.
\end{align*}
\end{dfn}
It is emphasized here that the weighted adjacency matrix $A(\mathcal{N})$ is a min-plus matrix.
Conversely, for any matrix $A\in\mathbb{R}_{\min}^{n\times n}$, 
there exists a network whose weighted adjacency matrix coincides with $A$.
We hereinafter denote such a network by $\mathcal{N}(A)$.
\par
If the vertex indices $i(0),i(1),\dots,i(s)$ are different from each other, 
and edges $e_{i(0),i(1)}, e_{i(1),i(2)},\dots,e_{i(s-1),i(s)}$ exist, 
then $P=(v_{i(0)},v_{i(1)},\dots,v_{i(s)})$ is a path on the network $\mathcal{N}$.
For the path $P$, the length $\ell (P)$ denotes the edge number $s$, 
and the weight sum $\omega (P)$ designates the sum of the edge weights:
\begin{align*}
\omega(P)=\sum_{k=0}^{s-1}w((v_{i(k)},v_{i(k+1)}))=\sum_{k=0}^{s-1}a_{i(k)i(k+1)}
=\bigotimes_{k=0}^{s-1}a_{i(k)i(k+1)}. 
\end{align*}
Moreover, the path $P$ with $i(0)=i(s)$ is just a circuit, 
and its length and weight are calculated in the same manner as those in path $P$.
The following definition describes the average weight of the circuit $C$. 
\begin{dfn}
For the circuit $C$, the average weight $\text{ave}(C)$ is given by 
\begin{align*}
\text{ave}(C)=\dfrac{\omega(C)}{\ell(C)}.
\end{align*}
\end{dfn}
The eigenvalues and eigenvectors of matrices play important roles in both linear algebra and min-plus algebra.
The following definition determines the eigenvalues and eigenvectors of the min-plus matrix.
\begin{dfn}
For the min-plus matrix $A\in\mathbb{R}_{\min}^{n\times n}$, if there exist $\lambda\in\mathbb{R}_{\min}$ and 
$\bm{x}\in\mathbb{R}_{\min}^n\!\setminus\!\{(\varepsilon,\varepsilon,\dots,\varepsilon)^{\top}\}$ satisfying 
\begin{align*}
A\otimes\bm{x}=\lambda\otimes\bm{x},
\end{align*}
then $\lambda$ and $\bm{x}$ are an eigenvalue and its corresponding eigenvector.
\end{dfn}
The eigenvalues of the min-plus matrix were shown in Baccelli et al. \cite{BCOQ} and Gondran-Minoux \cite{GM} 
to have interesting relationships with circuits in the network.
\begin{thm}[Baccelli et al \cite{BCOQ} and Gondran-Minoux \cite{GM}]
If the min-plus matrix $A\in\mathbb{R}_{\min}^{n\times n}$ has an eigenvalue $\lambda\not=\varepsilon$, 
there exists a circuit in the network $\mathcal{N}(A)$ whose average weight is equal to $\lambda$.
\end{thm}
\begin{thm}[Baccelli et al \cite{BCOQ} and Gondran-Minoux \cite{GM}]\label{prp;mineigen}
The minimum of the average weights of circuits in the network $\mathcal{N}(A)$ 
coincides with the minimum eigenvalue of the min-plus matrix $A\in\mathbb{R}_{\min}^{n\times n}$.
\end{thm}
In particular, Theorem \ref{prp;mineigen} suggests that we can algebraically compute 
the minimum of average weights of circuits in $\mathcal{N}(A)$ without grasping pictorial situations.
%
\section{Factorization of min-plus polynomials}
%
In this section, we briefly review Maclagan-Sturmfels \cite{MS} with regards to linear factorization over the min-plus algebra, 
and then describe a preconditioning algorithm in linear factorizations which will be helpful in later sections.
\par
\begin{figure}[htbp]
\centering
\includegraphics[width=80mm]{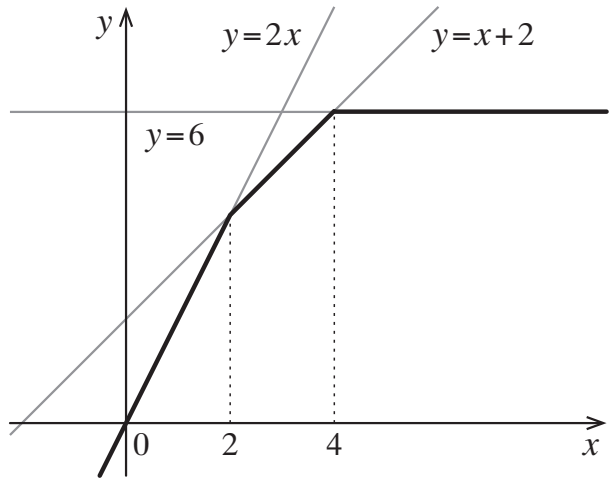}
\caption{The functional graph of $x^2\oplus 2\otimes x\oplus 6$.}
\label{graph1}
\end{figure}

We now consider the so-called min-plus polynomial of degree $n$ with respect to $x$,
\begin{align*}
p(x)=x^{n}\oplus c_1\otimes x^{n-1}\oplus\cdots\oplus c_{n-1}\otimes x\oplus c_n.  
\end{align*}
where $x^k:=\underbrace{x\otimes x\otimes\cdots\otimes x}_{k\text{ times}}=kx$ 
and $c_1,c_2,\dots,c_n\in\mathbb{R}_{\min} $ are the coefficients.
The following proposition gives the necessary and sufficient condition for factorizing the min-plus polynomial $p(x)$ 
into linear factors as $p(x)=(x\oplus c_1)\otimes [x\oplus (c_2-c_1)]\otimes\cdots\otimes [x\oplus (c_n-c_{n-1})]$.
\begin{prp}[Maclagan-Sturmfels \cite{MS}]\label{prp1}
The min-plus polynomial $p(x)$ can be completely 
factorized into linear factors if, and only if, the coefficients 
$a_0,a_1,\dots,a_{n-1}\in\mathbb{R}_{\min}$ satisfy the following inequality,
\begin{align*}
c_1\le c_2-c_1\le\cdots\le c_n-c_{n-1}.
\label{eq-lin-fact}
\end{align*}
\end{prp}
Regarding $p(x)$ as the min-plus function with respect to $x$, we see that $p(x)$ is piecewise linear.
This is because 
\begin{align*}
p(x)=\min\{nx,c_1+(n-1),\dots,c_{n-1}+x,c_n\}.
\end{align*}
Thus, the functional graph consists of line segments and rays.
Figure \ref{graph1} shows an example of the functional graph of $x^2\oplus 2\otimes x\oplus 6$.
\par
The piecewise linearity implies that $p(x)$ has a finite number of break points.
Roots of the min-plus function $p(x)$ coincide with the values of the $x$-coordinates of break points.
From Figure \ref{graph1}, we can thus factorize $p(x)$ as $p(x)=(x\oplus 2)\otimes (x\oplus 4)$.
It is emphasized here that two distinct min-plus polynomials, 
$p(x)$ and $p'(x)$, are sometimes factorized using common linear factors, which differs from over linear algebra.
If the linear factorizations of $p(x)$ and $p'(x)$ are the same, 
then we recognize that $p(x)$ is equivalent to $p'(x)$.
To distinguish $p(x)=q(x)$, namely, $p(x)$ is completely equal to $p'(x)$, 
we express $p(x)\equiv p'(x)$ if $p(x)$ is equivalent to $p'(x)$.
\par
We later need to find equivalent min-plus polynomials 
to observe the characteristic polynomials of min-plus matrices.
Although it is not so difficult to derive equivalent min-plus polynomials, 
we show how to reduce them to equivalent ones that can be directly factorized into linear factors.
Such algorithms, to the best of our knowledge, have not previously been presented.
Therefore, we describe an algorithm for constructing 
an equivalent polynomial that can be factorized into linear factors. 
\begin{alg}\label{alg1}
Constructing $p'(x)= x^n\oplus c'_1\otimes x^{n-1}\oplus\cdots\oplus c'_{n-1}\otimes x\oplus c'_n$
which is equivalent to $p(x)=x^n\oplus c_1\otimes x^{n-1}\oplus\cdots\oplus c_{n-1}\otimes x\oplus c_n$, 
namely, $p(x)\equiv p'(x)$.\\
{\bf Input}:
The coefficients $c_1,c_2,\dots,c_n$ in the min-plus polynomial $p(x)$.\\
{\bf Output}: 
The coefficients $c'_1,c'_2,\dots,c'_n$ in the equivalent min-plus polynomial $p'(x)$.\\
01: 
Set $c_1=0$ and $i:=0$.\\
02: 
Set $c_j=\varepsilon$ if $p(x)$ does not involve $x^{n-j}$.\\
03:
Compute $T_k:=(c_k -c_i)/(k-i)$ for $k=i+1,i+2,\dots,n$.\\
04:
Find integer $m$ such that $T_m=\displaystyle\min_{k=i+1,i+2,\dots,n} T_k$.\\
05:
Compute $c'_{i+1},c'_{i+2},\dots,c'_{m}$ as 
\begin{align*}
c'_\ell=\left\{\begin{array}{l}
c_i+(\ell-i)\dfrac{c_m-c_i}{m-i},\quad\ell=i+1,i+2,\dots,m-1,\\
c_m,\quad\ell=m.
\end{array}\right.
\end{align*}
06:
Overwrite $i:=m$.\\
07:
Set $c'_n:=c_n$ if $i=n$. Otherwise, go back to line 03.
\end{alg}
%
\section{All roots of characteristic polynomials of min-plus matrices}
Characteristic polynomials of matrices over linear algebra have roots which are just the eigenvalues.
However, to the best of our knowledge, the characteristic polynomials of min-plus matrices have not yet been strictly defined.
Min-plus characteristic polynomials can be, for example, given using the tropical determinant.
Such characteristic polynomials have minimum roots which coincide with minimum eigenvalues 
and the minimums of average weights in the corresponding networks.
The literature has not discussed whether the other roots are eigenvalues or not, 
nor whether they are meaningful features or not in the network.
In this section, we thus clarify the relationship between the 2nd, 3rd,\dots, 
minimum roots and the average weights of circuits in a special network.
\par
We first review characteristic polynomials of min-plus matrices using the tropical determinant \cite{MS}.
For the min-plus matrix $A=(a_{ij})\in\mathbb{R}_{\min}^{n\times n}$, 
the tropical determinant, denoted $\text{tropdet}(A)$, is defined by:
\begin{align*}
\text{tropdet}(A)=\bigoplus_{\sigma\in S_n}a_{1\sigma(1)}\otimes a_{2\sigma(2)}\otimes\cdots\otimes a_{n\sigma(n)},
\end{align*}
where $S_n$ is the symmetric group of permutations of $\{1,2,\dots,n\}$.
The following definition then determines the characteristic polynomial of $A$.
\begin{dfn}\label{trop}
For the min-plus matrix $A\in\mathbb{R}_{\min}^{n\times n}$, the characteristic polynomial $g_A(x)$ is given by 
\begin{align*}
g_A(x)=\text{tropdet}(A\oplus x\otimes I), 
\end{align*}
where $I$ is the $n$-by-$n$ identity matrix 
whose $(i,j)$ entries are $0$ if $i=j$, or $\varepsilon$ otherwise.
\end{dfn}
\par
To distinguish the distinct circuits in the network $\mathcal{N}(A)$, 
that are associated with the min-plus matrix $A\in\mathbb{R}_{\min}^{n\times n}$, 
we hereinafter use the notation $C(\ell_i,p_i)$ 
as the circuit of length $\ell_i$ and with the average weight $p_i$.
Moreover, we prepare a set of circuits with a length sum of $\tilde{\ell}_i$ in the network $\mathcal{N}(A)$.
Here, we regard the extended circuit of length $\tilde{\ell}_i$, and designate it 
as $\tilde{C}(\tilde{\ell}_i,\tilde{p}_i)$ where $\tilde{p}_i$ is the average weight.
Of course, simple circuits are members of extended circuits, 
and the weight sum of $\tilde{C}(\tilde{\ell}_i,\tilde{p}_i)$ is $\tilde{\ell}_i\tilde{p}_i$ for each $i$.
According to Maclagan-Sturmfels \cite{MS},  we can easily derive a proposition 
concerning the relationships between coefficients of the characteristic polynomial 
and the weight sums of extended circuits in the network.
\begin{prp} \label{prp-eigen}
For the min-plus matrix $A\in\mathbb{R}_{\min}^{n\times n}$, 
let us assume that the characteristic polynomial $g_A(x)$ is expanded as 
\begin{align*}
g_A(x)\equiv x^n\oplus c_1\otimes x^{n-1}\oplus\cdots\oplus c_{n-1}\otimes x\oplus c_n.
\end{align*}
Then, each coefficient $c_j$ coincides with the minimum of the weight sums of the extended circuits 
in the set of the separated and extended circuits 
$\mathcal{C}_j:=\{\tilde{C}(\tilde{\ell}_i,\cdot)\mid\tilde{\ell}_i=j\}$ 
in the network $\mathcal{N}(A)$ that are associated with $A$.
\end{prp}
Now, we consider the case where $k$ separate circuits $C(\ell_1,p_1),C(\ell_2,p_2),\dots,$ $C(\ell_k,p_k)$ 
existin the network $\mathcal{N}$.
Strictly speaking, $C(\ell_1,p_1),C(\ell_2,p_2),\dots,$ $C(\ell_k,p_k)$ 
are distinct to each other and every vertex belongs to at most one circuit in the network $\mathcal{N}$.
Without loss of generality, we may assume that $p_1\le p_2\le\cdots\le p_k$.
Moreover, we recognize that the extended circuit  $\tilde{C}(\tilde{\ell}_i,\tilde{p}_i)$ is homogeneous 
if all simple circuits in $\tilde{C}(\tilde{\ell}_i,\tilde{p}_i)$ have the same average weight $\tilde{p}_i$.
We then see that, in the case where $C(\ell_1,p_1),$ $C(\ell_2,p_2),\dots,C(\ell_k,p_k)$ are separate circuits, 
$j$ homogeneous extended circuits $\tilde{C}(\tilde{\ell}_1,\tilde{p}_1),\tilde{C}(\tilde{\ell}_2,\tilde{p}_2),\dots,$ 
$\tilde{C}(\tilde{\ell}_j,\tilde{p}_j)$ exist where $\tilde{p}_1< \tilde{p}_2<\cdots<\tilde{p}_j$ 
and $j\le k$ in the network $\mathcal{N}$.
This is key role to deriving the following two main theorems in this section.
\begin{thm}\label{thm1}
Let us assume that all circuits are separated in the network $\mathcal{N}(A)$ 
associated with the min-plus matrix $A\in\mathbb{R}_{\min}^{n\times n}$.
Then the characteristic polynomial $g_A(x)$ can be factorized into linear factors of the form 
$g_A(x)\equiv (x\oplus\tilde{p}_1)^{\tilde{\ell}_1}\otimes (x\oplus\tilde{p}_2)^{\tilde{\ell}_2}\otimes\cdots
\otimes (x\oplus\tilde{p}_k)^{\tilde{\ell}_k}\otimes x^r$, 
where $r:=n-(\tilde{\ell}_1+\tilde{\ell}_2+\cdots+\tilde{\ell}_k)$.
\end{thm}
\begin{prf}
Without loss of generality, we may assume $\tilde{p}_1<\tilde{p}_2<\cdots<\tilde{p}_k$.
We first prove that $\tilde{p}_1,\tilde{p}_2,\dots,\tilde{p}_k$ are roots of 
$g_A(x)=x^n\oplus c_1\otimes x^{n-1}\oplus\cdots\oplus c_{n-1}\otimes x\oplus c_n$.
It is obvious that the leading term $x^n$ becomes $n\tilde{p}_1$ at $x=\tilde{p}_1$.
From Proposition \ref{prp-eigen}, the coefficient $c_{\tilde{\ell}_1}$ is equal to 
the minimum of weight sums of the extended circuits in the set $\mathcal{C}_{\tilde{\ell}_1}$.
Since $\tilde{p}_1$ is the minimum average weight, $c_{\tilde{\ell}_1}=\tilde{\ell}_1\tilde{p}_1$.
Thus, we can simplify the term $c_{\tilde{\ell}_1}\otimes x^{n-\tilde{\ell}_1}$ as 
$\tilde{\ell}_1\tilde{p}_1+(n-\tilde{\ell}_1)\tilde{p}_1=n\tilde{p}_1$ at $x=\tilde{p}_1$.
Similarly, for all $i\not=\tilde{\ell}_1$, $c_{i}\otimes x^{n-i}=c_i+(n-i)\tilde{p}_1$ at $x=\tilde{p}_1$.
If $c_i+(n-i)\tilde{p}_1<n\tilde{p}_1$, namely, $c_i/i<p_1$, 
then the homogeneous extended circuit $\tilde{C}(i,\tilde{p}_0)$ exists where $\tilde{p}_0<\tilde{p}_1$.
This contradicts the assumption that $\tilde{p}_1$ is the minimum average weight.
Thus, we conclude that $x=\tilde{p}_1$ is a root of $g_A(x)$.
Moreover, we can easily derive $c_{\tilde{\ell}_1+\tilde{\ell}_2}\otimes x^{n-\tilde{\ell}_1-\tilde{\ell}_2}
=\tilde{\ell}_1\tilde{p}_1+(n-\tilde{\ell}_1)\tilde{p}_2$ at $x=\tilde{p}_2$.
This is because Proposition \ref{prp-eigen} immediately leads to 
$c_{\tilde{\ell}_1+\tilde{\ell}_2}=\tilde{\ell}_1\tilde{p}_1+\tilde{\ell}_2\tilde{p}_2$.
Simultaneously, we can observe that $c_{\tilde{\ell}_1}+x^{n-\tilde{\ell}_1}
=\tilde{\ell}_1\tilde{p}_1+(n-\tilde{\ell}_1)\tilde{p}_2$.
Thus, to prove that $\tilde{p}_2$ is a root of $g_A(x)$, it is necessary to show that, for all 
$i\not=\tilde{\ell}_1,\tilde{\ell}_1+\tilde{\ell}_2$, $c_i\otimes x^{n-i}\ge\tilde{p}_1\tilde{\ell}_1
+ (n-\tilde{\ell}_1) \tilde{p}_2$ at $x=\tilde{p}_2$, namely, 
$c_i-i\tilde{p}_2 \ge\tilde{\ell}_1(\tilde{p}_1-\tilde{p}_2)$.
Recalling here that $\tilde{p}_1<\tilde{p}_2$, we see that $c_i-i\tilde{p}_2<0$, namely, 
$c_i/i<\tilde{p}_2$ if $c_i-i\tilde{p}_2<\tilde{\ell}_1(\tilde{p}_1-\tilde{p}_2)$.
This implies that $c_i/i=\tilde{p}_1$ for $i\not=\tilde{\ell}_1$, 
but $c_i/i\not=\tilde{p}_1$ for $i\not=\tilde{\ell}_1$.
Therefore, we recognize that $\tilde{p}_2$ is also a root of $g_A(x)$.
Along the same lines, we observe that, for $m=2,3,\dots,k$, only two terms: 
$\tilde{\ell}_1\tilde{p}_1\otimes\tilde{\ell}_2\tilde{p}_2\otimes\cdots\otimes\tilde{\ell}_{m-1}\tilde{p}_{m-1}
\otimes x^{n-\tilde{\ell}_1-\tilde{\ell}_2-\cdots-\tilde{\ell}_{m-1}}$
and $\tilde{\ell}_1\tilde{p}_1\otimes\tilde{\ell}_2\tilde{p}_2\otimes\cdots\otimes\tilde{\ell}_m\tilde{p}_m
\otimes x^{n-\tilde{\ell}_1-\tilde{\ell}_2-\cdots-\tilde{\ell}_{m}}$ become both 
$\tilde{\ell}_1\tilde{p}_1+\cdots+\tilde{\ell}_{m-1}\tilde{p}_{m-1}
+(n-\tilde{\ell}_1-\tilde{\ell}_2-\cdots-\tilde{\ell}_{m-1})$ at $x=\tilde{p}_m$,
and are the minimum among all terms in $g_A(x)$.
This suggests that $x=\tilde{p}_m$ is a root of $g_A(x)$.
\par
Next, we examine the linear factorization of $g_A(x)$.
We can update $c_1,c_2,\dots,c_{\tilde{\ell}_1}$ as 
$\tilde{p}_1,2\tilde{p}_1,\dots,\tilde{\ell}_1\tilde{p}_1$, respectively, using Algorithm \ref{alg1}.
We then see that $x^n,c_1\otimes x^{n-1},\dots,c_{\tilde{\ell}_1}\otimes x^{n-\tilde{\ell}_1}$ 
are equal to each other at $x=\tilde{p}_1$.
Similarly, Algorithm \ref{alg1} updates $c_{\tilde{\ell}_1+1},c_{\tilde{\ell}_1+2},\dots,
c_{\tilde{\ell}_1+\tilde{\ell}_2}$ as $\tilde{\ell}_1\tilde{p}_1+\tilde{p}_2,\tilde{\ell}_1\tilde{p}_1
+2\tilde{p}_2,\dots,\tilde{\ell}_1\tilde{p}_1+\tilde{\ell}_2\tilde{p}_2$, then, it holds that 
$c_{\tilde{\ell}_1}\otimes x^{n-\tilde{\ell}_1}=c_{\tilde{\ell}_1+1}\otimes x^{n-\tilde{\ell}_1-1}=\cdots
=c_{\tilde{\ell}_1+\tilde{\ell}_2}\otimes x^{n-\tilde{\ell}_1-\tilde{\ell}_2}$ at $x=\tilde{p}_2$. 
Applying Algorithm \ref{alg1} repeatedly, we see that 
\begin{align*}
c_{i+1}-c_i=\left\{\begin{array}{ll}
\tilde{p}_1, & i=0,1,\dots,\tilde{\ell}_1-1, \\
\tilde{p}_2, & i=\tilde{\ell}_1,\tilde{\ell}_1+1,\dots,\tilde{\ell}_1+\tilde{\ell}_2-1,\\
 & \qquad \vdots \\
\tilde{p}_k, & i=\tilde{\ell}_1+\tilde{\ell}_2+\cdots+\tilde{\ell}_{k-1}, 
\dots,\tilde{\ell}_1+\tilde{\ell}_2+\cdots+\tilde{\ell}_{k-1}+\tilde{\ell}_k-1,
\end{array}\right.
\end{align*}
where $c_0=0$. 
If $r=n-(\tilde{\ell}_1+\tilde{\ell}_2+\cdots+\tilde{\ell}_k)=0$, 
then it immediately follows from Proposition \ref{prp1} that 
$g_A(x)=(x\oplus\tilde{p}_1)^{\tilde{\ell}_1}\otimes (x\oplus\tilde{p}_2)^{\tilde{\ell}_2}
\otimes\cdots\otimes (x\oplus\tilde{p}_k)^{\tilde{\ell}_k}$. 
If $r>0$, then there is no extended circuits greater in length than 
$n-r$ in the network $\mathcal{N}(A)$.
This is because there exist $r$ vertices that do not belong to any circuits.
Thus, we can overwrite the coefficients $c_{n-r+1},c_{n-r+2},\dots,c_{n-1}$ and 
the constant term $c_n$ with 0.
Therefore, we have $g_A(x)=(x\oplus\tilde{p}_1)^{\tilde{\ell}_1}\otimes 
(x\oplus\tilde{p}_2)^{\tilde{\ell}_2}\otimes\cdots\otimes (x\oplus\tilde{p}_k)^{\tilde{\ell}_k}\otimes x^r $. \qed
\end{prf}
\begin{thm}\label{thm2}
For the min-plus matrix $A\in\mathbb{R}_{\min}^{n\times n}$, assume that all circuits are separated 
in the network $\mathcal{N}(A)$ associated with $A$.
If the characteristic polynomial $g_A(x)$ can be factorized into linear factors of the form 
$g_A(x)\equiv (x\oplus\tilde{p}_1)^{\tilde{\ell}_1}\otimes (x\oplus\tilde{p}_2)^{\tilde{\ell}_2}\otimes\cdots
\otimes (x\oplus\tilde{p}_k)^{\tilde{\ell}_k}\otimes x^r$, then there exist homogeneous extended circuits 
$\tilde{C}(\tilde{\ell}_1,\tilde{p}_1),\tilde{C}(\tilde{\ell}_2,\tilde{p}_2),
\dots,\tilde{C}(\tilde{\ell}_k,\tilde{p}_k)$.
\end{thm}
\begin{prf}
Similarly to prove Theorem \ref{thm1}, 
assume that $\tilde{p}_1<\tilde{p}_2<\cdots<\tilde{p}_k$.
Here, we focus on the case $r=0$.
Going over the proof of Theorem \ref{thm1}, we see that $g_A(x)$ is equivalent to
\begin{align*}
\hat{g}_A(x)&=x^n\oplus\tilde{\ell}_1\tilde{p}_1\otimes x^{n-\tilde{\ell}_1}\oplus 
(\tilde{\ell}_1\tilde{p}_1\otimes\tilde{\ell}_2\tilde{p}_2)\otimes 
x^{n-\tilde{\ell}_1-\tilde{\ell}_2}\oplus\cdots\\
&\quad\oplus (\tilde{\ell}_1\tilde{p}_1\otimes\cdots\otimes\tilde{\ell}_{k-1}\tilde{p}_{k-1})
\otimes x^{\tilde{\ell}_k}\oplus (\tilde{\ell}_1\tilde{p}_1\otimes\cdots\otimes\tilde{\ell}_k\tilde{p}_k)
\end{align*} 
The coefficients $\tilde{\ell}_1\tilde{p}_1,\tilde{\ell}_1\tilde{p}_1\otimes\tilde{\ell}_2\tilde{p}_2,
\dots,\tilde{\ell}_1\tilde{p}_1\otimes\tilde{\ell}_2\tilde{p}_2\otimes\cdots\otimes\tilde{\ell}_{k-1}\tilde{p}_{k-1}$
and the constant term $\tilde{\ell}_1\tilde{p}_1\otimes\tilde{\ell}_2\tilde{p}_2
\otimes\cdots\otimes\tilde{\ell}_k\tilde{p}_k$ imply that the network $\mathcal{N}(A)$ includes 
the homogeneous extended circuits $\tilde{C}(\tilde{\ell}_1,\tilde{p}_1),
\tilde{C}(\tilde{\ell}_2,\tilde{p}_2),\dots,\tilde{C}(\tilde{\ell}_k,\tilde{p}_k)$. \qed
\end{prf}
From Theorems \ref{thm1} and \ref{thm2}, we can conclude that the $2$nd, $3$rd, \dots $k$th 
minimum roots of the characteristic polynomial $g_{A}(x)$ are equal to the average weights 
$\tilde{p}_2,\tilde{p}_3,\dots,\tilde{p}_k$, respectively, if, and only if, 
the circuits $C(\ell_1,\tilde{p}_1),$ $C(\ell_2,\tilde{p}_2),\dots,C(\ell_k,\tilde{p}_k)$ are all separated.
%
\section{New characteristic polynomials}
In this section, we propose new characteristic polynomials of min-plus matrices 
by imagining the analogue of the Faddeev-LeVerrier algorithm \cite{Fad},
which is an algorithm for generating characteristic polynomials of matrices in linear algebra.
\par
In the Faddeev-LeVerrier algorithm, only the sums and products of scalars and matrices 
construct the characteristic polynomials of linear matrices.
In fact, for a linear matrix $A\in\mathbb{R}^{n\times n}$, the coefficients $c_1,c_2,\dots,c_n$ 
appearing in the characteristic polynomial $x^n+c_1x^{n-1}+\cdots+c_{n-1}x+c_n$ 
is recursively given as 
\begin{align*}
& c_1=-\text{Tr}(A),\\
& c_2=-\frac{1}{2}\text{Tr}(A^2+c_1A),\\
& \quad\vdots\\
& c_n=-\frac{1}{n}\text{Tr}(A^n+c_1A^{n-1}+\cdots+c_{n-1}A).
\end{align*}
Thus, we can derive new characteristic polynomials of min-plus matrices based on this method.
\begin{dfn}\label{frame}
For the min-plus matrix $A\in\mathbb{R}_{\min}^{n\times n}$, the characteristic polynomial $\hat{g}_A(x)$
\begin{align*}
\hat{g}_A(x)=x^n\oplus c_1\otimes x^{n-1}\oplus\cdots\oplus c_{n-1}\otimes x\oplus c_n, 
\end{align*}
is recursively given as 
\begin{align*}
& c_1=\text{Tr}(A),\\
& c_2=\text{Tr}(A^2\oplus c_1\otimes A),\\
& \quad\vdots\\
& c_n=\text{Tr}(A^n\oplus c_1\otimes A^{n-1}\oplus\cdots\oplus c_{n-1}\otimes A),
\end{align*}
where $A^k=A^{k-1}\otimes A$ for $k=2,3,\dots,n$.
\end{dfn}
It is remarkable that the new characteristic polynomial $\hat{g}_A(x)$ usually differs from 
the already known characteristic polynomial $g_A(x)$.
The following theorem gives the relationship between the minimum root of 
the characteristic polynomial $\hat{g}_A(x)$ and
the eigenvalue of the min-plus matrix $A\in\mathbb{R}_{\min}^{n\times n}$.
\begin{thm}\label{thm3}
For the min-plus matrix $A\in\mathbb{R}_{\min}^{n\times n}$, 
the minimum root of the characteristic polynomial $\hat{g}_A(x)$ is equal to the eigenvalue of $A$.
\end{thm}
\begin{prf}
With the help of Proposition \ref{prp;mineigen}, 
we may prove that the minimum root, denoted $p_{\min}$, 
is just the minimum of average weights of circuits in the network $\mathcal{N}(A)$.
Regarding $\hat{g}_A(x)$ as the function with respect to $x$, 
we recall that $p_{\min}$ coincides with the minimum of the $x$-coordinates 
of breakpoints on the corresponding $xy$ functional graph.
It is worth noting here that the breakpoints with the minimum $x$-coordinate are the intersection 
of two lines $y=nx$ and $y=c_i(n-i)x$ for some $i$.
Thus, we derive $p_{\min}=c_i/i$.
\par
It remains to proven that $c_i/i$ becomes the minimum of the average weights 
of circuits in the network $\mathcal{N}(A)$.
Obviously, the coefficient $c_1= \text{Tr}(A)$ is equal to the minimum of the weight sums 
of circuits in the set $\mathcal{C}_1$.
Taking into account that the diagonals of $A^2$ and $c_1\otimes A$ are the weight sums 
of all extended circuits in $\mathcal{C}_2$, we see that $c_2=\text{Tr}(A^2\oplus c_1\otimes A)$ 
expresses the minimum of the weight sums of all extended circuits in $\mathcal{C}_2$.
Similarly, $c_i$ signifies the minimum of the weight sums 
of all extended circuits in $\mathcal{C}_i$.
Thus, $p_{\min}=\min_ic_i/i$ is equal to the minimum of the average weights 
of all extended circuits in $\mathcal{C}_i$.
Simultaneously, we see that the average weights of all extended circuits 
in the network $\mathcal{N}(A)$ are equal to or larger than $p_{\min}$.
Moreover, if the minimum of the average weights of extended circuits in $\mathcal{C}_{\tilde{\ell}_i}$ 
is $p_{\min}$, then $\tilde{C}(\tilde{\ell}_i,\tilde{p}_i)$ is a simple circuit.
This is because, if $\tilde{C}(\tilde{\ell}_i,\tilde{p}_i)$ is not a simple circuit, namely, 
$\tilde{C}(\tilde{\ell}_i,\tilde{p}_i)= \{C(\ell_1,p_1),C(\ell_2,p_2),\dots,C(\ell_k,p_j)\}$, 
then the average weight $\min_{i=1,2,\dots,j}\{\text{ave}(C(\ell_j,p_j))\}$ is smaller than $p_{\min}$.
Therefore, we conclude that $p_{\min}$ becomes the minimum of the average weights 
of the circuits in the network $\mathcal{N}(A)$. \qed
\end{prf}
Although two characteristic polynomials, $g_A(x)$ and $\hat{g}_A(x)$ are essentially distinct, 
they are equivalent to each other in a special case.
\begin{crl}\label{crl1}
If all circuits are simple and separated in the network $\mathcal{N}(A)$, 
then two characteristic polynomials $g_A(x)$ and $\hat{g}_A(x)$ satisfy $g_A(x)\equiv\hat{g}_A(x)$.
\end{crl}
Corollary \ref{crl1} can be provided via the proofs of Theorems \ref{thm1}, \ref{thm2} and \ref{thm3}.
We here give an example to illustrate the difference 
between two characteristic polynomials $g_A(x)$ and $\hat{g}_A(x)$.
For the min-plus matrix 
\begin{align*}
A=\begin{pmatrix}
\varepsilon & \varepsilon & 2 & \varepsilon & \varepsilon & \varepsilon & \varepsilon \\
3 & \varepsilon & \varepsilon & 2 & \varepsilon & \varepsilon & \varepsilon \\
\varepsilon & 1 & 3 & 9 & 1 & \varepsilon & \varepsilon \\
\varepsilon & 6 & \varepsilon & \varepsilon & \varepsilon & 2 & \varepsilon \\
\varepsilon & \varepsilon & \varepsilon & \varepsilon & \varepsilon & 2 & 1 \\
\varepsilon & \varepsilon & \varepsilon & \varepsilon & \varepsilon & \varepsilon & 1 \\
\varepsilon & \varepsilon & \varepsilon & \varepsilon & \varepsilon & \varepsilon & \varepsilon
\end{pmatrix},
\end{align*}
we obtain two characteristic polynomials 
\begin{align*}
& g_A(x)=x^7\oplus 3\otimes x^6\oplus 8\otimes x^5\oplus 6\otimes x^4\oplus 20\otimes x^3,\\
& \hat{g}_A(x)=x^7\oplus 3\otimes x^6\oplus 6\otimes x^5\oplus 6\otimes x^4\oplus 9\otimes x^3
\oplus 12\otimes x^2\oplus 12\otimes x\oplus 15.
\end{align*}
Using Algorithm \ref{alg1}, 
we can factorize $g_A(x)$ and $\hat{g}_A(x)$ as 
\begin{align*}
& g_A(x)\equiv (x\oplus 2)^3\otimes (x\oplus 14)\otimes x^3,\\
& \hat{g}_A(x)\equiv (x\oplus 2)^6\otimes (x\oplus 3).
\end{align*}
As shown in Theorems \ref{thm1}, \ref{thm2} and \ref{thm3}, 
the minimum roots of $g_A(x)$ and $\hat{g}_A(x)$ are certainly both the eigenvalue of $A$.
However, since the 2nd minimum roots of $g_A(x)$ and $\hat{g}_A(x)$ are 14 and 3, respectively, 
they are not equal to each other.
In actuality, limited to computing the eigenvalue of $A$, we can equivalently simplify $\hat{g}_A(x)$ as 
\begin{align*}
\breve{g}_A(x)=x^4\oplus 3\otimes x^3\oplus 6\otimes x^2\oplus 6\otimes x\oplus 9
\equiv (x\oplus 2)^3\otimes (x\oplus 3).
\end{align*}
In other words, it is not necessary to determine the coefficients 
$c_5=c_2\oplus 6=12,c_6=c_3\oplus 6=12,c_7=c_4\oplus 6=15$ to compute the eigenvalue.
Obviously, the linear factorization of $\breve{g}_A(x)$ is easier than that of $g_A(x)$.
New characteristic polynomials are thus expected to gain more advantage, as the matrix-size increases.
%
\section{Concluding remarks}
In this paper, we focused on all the roots of the already known characteristic polynomials of matrices, 
which are given from the links of vertices in networks on graphs, over min-plus algebra, 
and presented distinct new characteristic polynomials.
First, we briefly explained scalar and matrix arithmetic over min-plus algebra, 
the eigenvalues of min-plus matrices and the minimum average weights of circuits in networks, 
and the linear factorizations of min-plus polynomials. 
We then described a preconditioning algorithm for performing effective linear factorizations. 
Of course, the eigenvalues of min-plus matrices are the minimum roots of the already known characteristic polynomials. 
In other words, the minimum roots coincide with the minimum average weights of circuits in the corresponding networks.
Restricting the case to one where all circuits are completely separated in networks, 
we next showed that the $2$nd, $3$rd, $\dots$ minimum roots of the already known characteristic polynomials 
are just the $2$nd, $3$rd, $\dots$ minimum average weights, respectively. 
Finally, we propose new characteristic polynomials whose minimum roots are also the eigenvalues of min-plus matrices, 
and showed that they are equivalent to the already known characteristic polynomials 
if all circuits are completely separated in networks. 
We provided an example to verify the difference between the already known and proposed characteristic polynomials. 
The example simultaneously suggests that the proposed characteristic polynomials can be substantially reduced 
if the edge number is not large in the corresponding networks.
Thus, the proposed characteristic polynomials may be, so to speak, minimal polynomials.
Future work will focus on examining this aspect and designing reduction algorithms.\\[15pt]
{\bf\large Acknowledgements}\\[3pt]
This was partially supported by Grants-in-Aid for Scientific Research (C) No. 26400208 
from the Japan Society for the Promotion of Science.
\end{document}